\documentclass[a4paper,12pt ,fleqn]{article}
\usepackage{times}
\usepackage{amsthm, amssymb, amsmath}
\usepackage{epsfig}

\newcommand{\qkdf}{$q$-Kamp\'{e} de F\'{e}riet }

\newcommand{\hyp}{hypergeometric}

\newcommand{\sn}{Section\ }
\newcommand{\aph}{\alpha}
\newcommand{\btt}{\beta}
\newcommand{\gmm}{\gamma}

\newcommand{\lrdl}[3]{\left #1 {#3} \right #2}

\newcommand{\qhl}{$q$-hypergeometric}
\newcommand{\qsf}{$q$-shifted factorials }
\newcommand{\id}{identities}

\newcommand{\idy}{identity}

\newcommand{\rr}{Rogers-Ramanujan}
\newcommand{\slz}{Saalsch\"{u}tz}

\newcommand{\bib}{\bibitem }
\newcommand{\cireta}[3]{\textsc{#1}; {\it #2}, #3, \ \  to appear.}
\newcommand{\cire}[5]{\textsc{#1}; {\it #2}, #3\   (#4), #5.}
\newcommand{\cibo}[4]{\textsc{#1}; {\it #2}, #3\  (#4).}

\newcommand{\z}[2]{\genfrac{}{}{1pt}{0}{#1}{#2}}

\newcommand{\ez}[3]{#1_{r} = \z{#2}{#3}}
\newcommand{\bez}[2]{\btt_{(n, l)} = \z{#1}{#2}}
\newcommand{\bezk}[2]{\btt_{(n, l, k)} = \z{#1}{#2}}

\newcommand{\sfrq}[2]{\lrdl{(}{)}{#1;q}_{#2}}

\newcommand{\vwp}[3]{_{#1}W_{#2}\lrdl{(}{)}{#3}}

\newcommand{\kdf}[9]{ \Phi ^{#1}_{#2} \left[ 
\begin{array}{ccc}#3:& #4;& #5 \\ [8pt] #6:& #7;& #8 \end{array}
{;#9} \right ]}

\newcommand{\spdt}[1]{\prod^{\infty}_
{\substack{m = 1 \\ m \ncong #1}}
(1-q^{m})^{-1}}

\newcommand{\asmm}[2]{\sum^{#2}_{#1}}

\newcommand{\qalg}{$q$-analogue}

\newcommand{\qpss}{$q$-Pfaff-\slz\ sum ~\cite[ p. 237, (II. 12)]{gas}}
\newcommand{\sxpf}{terminating very-well-poised $_6\Phi_5$ 
sum ~\cite[ p. 238, (II. 21)]{gas}}
\newcommand{\qjeps}{Jackson's \qalg\ of Dougall's $_7 F_6$ 
sum ~\cite[ p. 238, (II. 22)]{gas}}

\newcommand{\bmltl}[2]{\begin{\multline} #1 \tag{#2}\end{multline}} 
\title{\bf BAILEY TYPE TRANSFORMS AND APPLICATIONS}
\author{\bf C.M. Joshi and Yashoverdhan Vyas\footnote{He is grateful 
to the \textbf{Council of Scientific and Industrial Research (CSIR), New 
Delhi, INDIA} for providing the Senior Research Fellowship.} }
\date{}
\begin{document}
\maketitle
\begin{abstract}
\noindent The aim of this paper is to establish new series 
transforms of Bailey type and to show that these Bailey type 
transforms work as efficiently as the classical one and give 
not only new $q$-hypergeometric identities, converting double 
or triple series into a very-well-poised $_{10}\Phi_9$ or
$_{12}\Phi_{11}$ series but can also be 
utilized to derive new double and triple series \rr\ type \id\ and 
corresponding infinite families of Rogers-Ramanujan 
type identities. 
\end{abstract}
\section{ Introduction} 
It is well-known that W. N. Bailey~\cite{bly1, bly2} established a series 
transform known as Bailey's transform and gave a mechanism to derive 
ordinary and \qhl\ \id\ and \rr\ type \id .
Using Bailey's transform, Slater  ~\cite{slt1, slt2} derived 130 \rr\ type \id .
 It was Andrews~\cite{and1, and2, and8} and~\cite{and7}, who exploited  the 
Bailey's transform in the form of Bailey pair and Bailey chains 
to show that all of the 130 identities given by Slater~\cite{slt1, slt2, slt3}
 can be embedded in infinite families of multiple series \rr \ 
type \id . \\[9pt]  \indent 
After Bailey~\cite{bly1, bly2} and notably after 
Andrews~\cite{and7},  a large number of Mathematicians have 
worked on Bailey's transform in the form of Bailey lemma and 
Bailey chain to  make applications in 
the theory of generalized \hyp \ series, number 
theory, partition theory, combinatorics, physics and 
computer-algebra (see~\cite{ag, akandbre, and3, 
and4, and5, and6, and9, and10, and11, and12, and13, 
andet1, andberch1, andberch2, abdbowman, andscisow, berchet1, 
berchet2, brs1, brs2, coskh, densingh, ext, foda, gas, jain1, jain2, 
jain3, lovjay1, paule1, paule2, paule3, acideka, acisow1, 
acisow2, acisow3, acisow4, acisow5,  sills1, sills2, 
sills3, sills4, sps, ubs, spird1, vj1, vj2, sow1, sow2, sow3, 
sow4, sow5, sow6}
 and~\cite{sow7} ) . \\[9pt]  \indent 
A survey of all the above mentioned papers reveals that the work done by 
Bailey~\cite{bly1, bly2} and Andrews~\cite{and7} may be considered as basic and 
fundamental in this field. 
At this point there arises one natural question, 
\begin{quotation}\noindent  {\it ``Can there be other series 
transforms of Bailey type, which  not only produce new 
ordinary and \qhl\ \id\ and \rr\ type \id\  but can also 
be iterated to provide multi sum versions of the \rr\ type 
\id ?" }
 \end{quotation} 
 This paper provides the answer in affirmative. \\[9pt]  \indent 
In this paper, first we establish two new series transform of
 Bailey type in \sn 2 as theorems 2.1 and 2.2 and then using 
these theorems we derive five new \qhl\ \id\ converting 
double or triple \qhl\ series into a very-well-poised $_{10}\Phi_9$ 
or $_{12}\Phi_{11}$ series in \sn 3. In \sn 4, we establish the
first Bailey type lemma and using it derive four new double 
series \rr\ type \id\ of modulo 7 and 5 and embed these in 
multiple series \rr\ type \id\ of modulo $4k+3$ and $4k+1$. 
In \sn 5, we establish second Bailey type lemma  and then 
using it derive two new triple series \rr\ type \id\ of 
modulo 9 and embed these  in multiple series \rr\ type \id\ 
of modulo $6s+3$. During the presentation of both lemmas in 
\sn 4 and 5, we establish only, the main results corresponding to 
Andrews~\cite{and7, and8} presentation of Bailey lemma.  \\[9pt]  \indent
It may be noted that the double and triple series \rr\ type \id\ presented 
in this paper do not arise as merely ``one or two level up" in the 
standard Bailey chain from well-known single sum identities.  
Furthermore, the original \rr\ \id\ (see~\cite{mac, ram, rog1} )
were also of modulus 5 but 
were single sum identities. A modulo 7 \idy\ with single sum is 
by Rogers~\cite[ Eq. 6; p. 331]{rog2} and as a one level up is by Andrews 
~\cite[ Eq. (1.8); p. 4083]{and2}. 
Andrews~\cite{and2} obtained his \idy\ of modulo 7 by particularizing 
the general multi sum odd moduli \id\ to double series case. 
After comparing our identities with the results of references 
~\cite{jain1, jain2, jain3, sills1, sills2, sills3, sills4, 
vj1, vj2} and so forth, we can say that
the \id\ derived in this paper are  new and most attractive and symmetric 
among all the  double and triple series \rr\ type \id\ derived 
previously. \\[9pt]  \indent 
The main tools in developing these transforms are two 
series rearrangements (2.5) and (2.10). Some of the instances of 
use of (2.5) and (2.10) may be found in the works of Burchnall 
and Chaundy~\cite{buch1}, Shankar~\cite{os1} 
and Srivastava and Manocha~\cite[p. 335]{hms1}.
The notation for double \qkdf\ in the eqn (3.1) and (3.3) is 
from~\cite[Eq. (282), (283); p. 349]{hms1} and of a triple series $\Phi^{(3)}$ in eqn (3.5) is simply 
a $q$-analogue of a particular case of the Srivastava's general triple series 
~\cite[Eq. (14), (15); p. 44]{hms1}. In view of notational 
difficulties mentioned in~\cite[pp. 270-274]{hms1}, the double non 
\qkdf\ series in eqns (3.2) and (3.4) have been written explicitly without 
using any notation. However, notation of $q$-analogue of Srivastava-Daoust's 
series~\cite[Eq. 284; p. 350]{hms1} may be used in these cases. Remaining definitions and 
notations are from~\cite{gas, and8} .
\section{Bailey type transforms}
\subsection{First Bailey type transform}
\subsubsection*{Theorem 2.1}
$$ \mbox{If  } \beta_{(n, l)} = \sum^{min.(n, l)}_{m=0} \alpha_m
 \ u_{n-m}\ u'_{l-m}  \ v_{n+m} \ v'_{l+m} 
\ t_{n-l} \ w_{l+n} \eqno (2.1) $$
$$ \mbox{ and }  \gamma_m = \sum^{\infty}_{n=m}
\sum^{\infty}_{l=m} \delta_n \ \delta'_l \
u_{n-m} \ u'_{l-m}  \ v_{n+m} \ v'_{l+m} 
\ t_{n-l} \ w_{l+n} \ \eqno (2.2) $$ then, 
subject to convergence conditions $$ \sum^ \infty_{m=0} \alpha_m 
\gamma_m = \sum^\infty_{n, l=0}\beta_{(n, l)} \delta_{n} 
\ \delta'_l , \eqno (2.3) $$
 where,  $\alpha_r, \delta_r, u_{r}, v_{r}, \delta'_r, u'_{r}, v'_{r}, 
t_{r}$ and $w_{r}$ are any functions of r only. 
 \subsubsection*{Proof :} 
Observe that $$ \sum^ \infty_{m=0} \alpha_m \ \gamma_m = \sum^ 
\infty_{m=0} \sum^{\infty}_{n=m}
\sum^{\infty}_{l=m} \aph_m\ \delta_n \ \delta'_l \
u_{n-m} \ u'_{l-m}  \ v_{n+m} \ v'_{l+m} 
\ t_{n-l} \ w_{l+n} . \eqno (2.4)$$ If this double series is 
convergent, then using  [30; p. 10, lemma 3], viz., $$ \sum^ 
{\infty}_{m=0} \sum^{\infty}_{n=0} \sum^{\infty}_{l=0} A(m, n, l) = 
\sum^{\infty}_{n=0} \sum^{\infty}_{l=0} \sum^{min. (n, l)}_{m=0}
 A(m, n-m, l-m) , \eqno (2.5)$$ in 
(2.4) after the replacement of $n$ by $n+m$ and of $l$ by $l+m$
, we get $$ 
\sum^ \infty_{m=0} \alpha_m \ \gamma_m =
\sum^{\infty}_{n=0} \sum^{\infty}_{l=0} \sum^{min. (n, l)}_{m=0}
 \aph_m\ \delta_n \ \delta'_l \
u_{n-m} \ u'_{l-m}  \ v_{n+m} \ v'_{l+m} 
\ t_{n-l} \ w_{l+n}$$ $$= \sum^\infty_{n, l=0}\beta_{(n, l)} \delta_{n} 
\ \delta'_l . $$
\subsection{Second Bailey type transform}
\subsubsection*{Theorem 2.2}
$$ \mbox{If  } \beta_{(n, l, k)} = 
\sum^{min.(n, l, k)}_{m=0} \alpha_m
 \ u_{n-m}\ u'_{l-m}\ u''_{k-m}  \ v_{n+m} \ v'_{l+m}\ v''_{k+m} 
 \eqno (2.6) $$
$$ \mbox{ and }  \gamma_m = \sum^{\infty}_{n=m}
\sum^{\infty}_{l=m} \sum^{\infty}_{k=m} 
\delta_n \ \delta'_l\ \delta''_k \
\ u_{n-m}\ u'_{l-m}\ u''_{k-m}  \ v_{n+m} \ v'_{l+m}\ v''_{k+m}
 \eqno (2.7) $$ then, 
subject to convergence conditions 
$$ \sum^ \infty_{m=0} \alpha_m 
\gamma_m = \sum^\infty_{n, l, k=0}\beta_{(n, l, k)} \delta_{n} 
\ \delta'_l\ \delta''_k \, \eqno (2.8) $$
 where,  $\alpha_r, \delta_r, u_{r}, v_{r}, \delta'_r, u'_{r}, v'_{r}, 
\delta''_r, u''_{r},$ and $ v''_{r} $ are any functions of r only. 
 \subsubsection*{Proof :} 
Observe that $$ \sum^ \infty_{m=0} \alpha_m \ \gamma_m = $$
$$\sum^ 
\infty_{m=0} \sum^{\infty}_{n=m}
\sum^{\infty}_{l=m}\sum^{\infty}_{k=m} 
 \aph_m\ \delta_n \ \delta'_l\ \delta''_k \
\ u_{n-m}\ u'_{l-m}\ u''_{k-m}  \ v_{n+m} \ v'_{l+m}\ v''_{k+m}
 .  \eqno (2.9)$$ If this double series is 
convergent, then using  [30; p. 10, lemma 3], viz., $$ \sum^ 
{\infty}_{m=0} \sum^{\infty}_{n, l, k = 0} A(m, n, l, k) = 
\sum^{\infty}_{n, l, k = 0} 
 \sum^{min. (n, l, k)}_{m=0}
 A(m, n-m, l-m, k-m) , \eqno (2.10)$$ in 
(2.8) after the replacement of $n, l$ and $k$ by $n+m, l+m$ 
and  $k+m$ respectively, we get $$ 
\sum^ \infty_{m=0} \alpha_m \ \gamma_m =$$
$$
\sum^{\infty}_{n=0} \sum^{\infty}_{l=0}\sum^{\infty}_{k=0}
 \sum^{min. (n, l, k)}_{m=0}
 \aph_m\ \delta_n \ \delta'_l\ \delta''_k \
\ u_{n-m}\ u'_{l-m}\ u''_{k-m}  \ v_{n+m} \ v'_{l+m}\ v''_{k+m}$$ 
$$
= \sum^\infty_{n, l, k=0}\beta_{(n, l, k)}\  \delta_{n} 
\ \delta'_l\ \delta''_k  . $$
Here, it may be noted that in
these two 
theorems  various new sequences of other combinations 
of involved summation indices may be introduced in their $\beta_{(n, l)}$ 
or $\beta_{(n, l, k)}$ and $\gamma_m $.
But we have presented these theorems with only those sequences, 
which have been used at least once in this paper.
 Further, it may be possible in many cases to decide 
the new expressions for $\alpha_r, \delta_r, u_{r}, v_{r}, \delta'_r, u'_{r}, v'_{r}, 
\delta''_r, u''_{r}, v''_{r}, t_{r}$ and $w_{r}$, 
which yield closed forms for $\beta_{(n, l)}$ or $\beta_{(n, l, k)}$ and 
$\gamma_m $, using one  or more of the known summation theorems. 
 Thus, many more new results may be discovered. But we have 
discussed certain cases which use the well-known classical 
summation theorems~\cite[ Appendix, II. 5, 6, 12, 21, 22]{gas} only. It may be emphasized that there may be other various appropriate cases to have the closed forms for 
$\beta_{(n, l)}$ or $\beta_{(n, l, k)}$ and $\gamma_m$ in theorems 2.1 and 2.2. \\[9pt]  \indent
Further, many series rearrangements, similar to (2.5) and (2.10), may also be utilized to discuss other Bailey type transforms  and their applications.
\section{ New $q$-\hyp \ \id } In the Bailey type 
transforms, discussed in the previous Section as theorems 2.1 and 2.2, we 
observed five expressions for $\alpha_r, \delta_r, u_{r}, v_{r}, \delta'_r, u'_{r}, v'_{r}, 
\delta''_r, u''_{r}, v''_{r}, t_{r}$ and $w_{r}$, 
which yield closed forms for $\beta_{(n, l)}$ or $\beta_{(n, l, k)}$ and 
$\gamma_m $ and lead to five new $q$-\hyp \ \id  . The \id \ 
are as follows:\\[9pt]   
\begin{equation*}
\kdf{1: 3; 2}{1: 2; 2}
{\z{qa}{d}} {b, c, q^{-M}} {b', c', q^{-N}}
{qa }{\z{qa}{d}, \z{b c q^{-M}}{a}}
{\z{qa}{d}, \z{b' c' q^{-N}}{a}}
{q, q}
\end{equation*}
$$
=\z{\sfrq{\z{qa}{b}, \z{qa}{c}}{M}}
{\sfrq{qa, \z{qa}{bc}}{M}}
\z{\sfrq{\z{qa}{b'}, \z{qa}{c'}}{N}}
{\sfrq{qa, \z{qa}{b'c'}}{N}}
\times\ \vwp{10}{9}{a; d, b, c, b', c', q^{-M}, q^{-N}; q, 
\z{a^{3}\ q^{3+M+N}}{b\ c\ b'\ c'\ d}} , 
\eqno (3.1)$$ \\[9pt]
$$
\asmm{n,l=0}{\infty} \z{\sfrq{\z{A}{a}}{n-l}}{\sfrq{\z{qa}{B}}{n-l}}
\z{\sfrq{\z{AB}{a}}{n+l}}{\sfrq{qa}{n+l}}
\z{\sfrq{A, q\sqrt{A}, -q\sqrt{A}, b, c, \z{Aaq^{1+M}}{b\ c}
,\z{qa}{B}, q^{-M}}{n}}
{\sfrq{\sqrt{A}, -\sqrt{A}, \z{qA}{b}, \z{qA}{c}, 
\z{b\ c\ q^{-M}}{a}, \z{AB}{a}, A\ q^{1+M}}{n}}
$$
$$
\times\ \z{\sfrq{B, q\sqrt{B}, -q\sqrt{B}, b', c', \z{Baq^{1+N}}{b'\ c'}
, \z{qa}{A}, q^{-N}}{l}}
{\sfrq{\sqrt{B}, -\sqrt{B}, \z{qB}{b'}, \z{qB}{c'}, 
\z{b'\ c'\ q^{-N}}{a}, \z{AB}{a}, B\ q^{1+N}}{l}}
\z{q^{n}\ \lrdl{(}{)}{\z{qAB}{a^{2}}}^l}{(q;q)_n \ (q;q)_l}
$$
$$
= \z{\sfrq{qA, \z{qa}{b},\z{qa}{c}, \z{qA}{bc} }{M}}
{\sfrq{qa, \z{qA}{b},\z{qA}{c}, \z{qa}{bc} }{M}} \cdot\
\z{\sfrq{qB, \z{qa}{b'},\z{qa}{c'}, \z{qB}{b'c'} }{N}}
{\sfrq{qa, \z{qB}{b'},\z{qB}{c'}, \z{qa}{b'c'} }{N}}
$$
$$
\times\ \vwp{12}{11}{a; \z{qa^{2}}{AB}, b, c, b', c', 
\z{Aaq^{1+M}}{b\ c}, \z{Baq^{1+N}}{b'\ c'}, q^{-M}, q^{-N};q,\ q} , 
\eqno(3.2)$$ \\[9pt]
$$
\kdf{2:2;8}{2:1;7}{b', \z{qa}{d}}{c', q^{-M}}
{\sqrt{qb'}, -\sqrt{qb'}, b, c, \z{ab'\ q^{1+N}}{b\ c}, 
\z{b'\ c'}{a}, \z{b'\ q^{-M}}{a}, q^{-N}}
{qa, \z{b'\ c'\ q^{-M}}{a}}{\z{qa}{d}}
{\sqrt{b'}, -\sqrt{b'}, \z{qb'}{b}, \z{qb'}{c}, \z{bcq^{-N}}{a}, 
b'q^{1+N}, \z{qa}{d}}{q,q}
$$
$$
=\z{\sfrq{\z{qa}{b'}, \z{qa}{c'}}{M}}{\sfrq{qa, \z{qa}{b'\ c'}}{M}}
\ \cdot\ 
\z{\sfrq{qb', \z{qa}{b}, \z{qa}{c}, \z{qb'}{bc}}{N}}
{\sfrq{qa, \z{qb'}{b}, \z{qb'}{c}, \z{qa}{bc}}{N}}$$
$$
\times\ \vwp{10}{9}{a; d, b, c, c', \z{q^{1+N}ab'}{bc}, q^{-M},
 q^{-N}; q, \z{q^{2+M}a^{2}}{b'\ c'}} \eqno(3.3)
$$\\[9pt]
 $$\hspace{-7cm}
 \asmm{n,l=0}{\infty} \z{\sfrq{\z{b'}{A}}{n-l}}{\sfrq{\z{qa}{A}}{n-l}}
\z{\sfrq{b'}{n+l}}{\sfrq{qa}{n+l}}
$$
$$\cdot\ \z{\sfrq{ q\sqrt{b'}, -q\sqrt{b'}, b, c, \z{b'\ q^{-M}}{a}, 
\z{b'\ c'}{a}, \z{b'\ a\ q^{1+N}}{b\ c}, \z{b'\ c'\ q^{-M}}{A}, 
\z{qa}{A}, q^{-N}}{n}}
{\sfrq{\sqrt{b'}, -\sqrt{b'}, \z{b'\ q^{-M}}{A}, 
\z{b'\ c'}{A}, \z{qb'}{b}, \z{qb'}{c}, \z{b'\ c'\ q^{-M}, 
}{a}, \z{bcq^{-N}}{a}, b'\ q^{1+N}}{n}}
 $$
 $$
\cdot\ \z{\sfrq{q \sqrt{A}, -q \sqrt{A},c', A, \z{Aa\ q^{1+M}}{b'\ c'}, 
q^{-M} }{l}}
{\sfrq{\sqrt{A}, -\sqrt{A}, \z{qA}{c'}, \z{b'\ c'\ q^{-M}}{a}, 
Aq^{1+M}}{l}} \cdot\ 
\z{q^{n}\ \lrdl{(}{)}{\z{q^{2+M}\ a^{2}}{b'\ c'}}^{l}}{(q;q)_n \ (q;q)_l} 
 $$
 $$
 =\z{\sfrq{qA, \z{qa}{b'}, \z{qa}{c'}, \z{qA}{b'\ c'}}{M}}
{\sfrq{qa, \z{qA}{b'}, \z{qA}{c'}, \z{qa}{b'\ c'}}{M}} \cdot\ 
\z{\sfrq{qb', \z{qa}{b}, \z{qa}{c}, \z{qb'}{bc}}{N}}
{\sfrq{qa, \z{qb'}{b}, \z{qb'}{c}, \z{qa}{bc}}{N}}
 $$
 $$
 \times\ \vwp{10}{9}{a; b, c, c', \z{a^{2}q^{1+M}}{b'\ c'},
 \z{ab'\ q^{1+N}}{bc}, q^{-M}, q^{-N};q,\ \z{qc'}{aq^{M}}}  
 \eqno(3.4)$$
\\ and 
$$
 \Phi^{(3)}\  \left[ 
\begin{array}{ccccccc}qa::& \hspace{3mm} :& \hspace{3mm} :& \hspace{3mm} :& 
b,c,q^{-N};& b',c',q^{-L};& b'',c'',q^{-K} \\ 
\\ \hspace{3mm} ::& qa :& qa :& qa:& \z{bcq^{-N}}{a};& 
\z{b'\ c'\ q^{-L}}{a};& \z{b''\ c''\ q^{-K}}{a} \end{array} 
{;q;q,q,q} 
\right]
$$ 
$$
=\z{\sfrq{\z{qa}{b}, \z{qa}{c}}{N}}{\sfrq{qa, \z{qa}{bc}}{N}}
\z{\sfrq{\z{qa}{b'}, \z{qa}{c'}}{L}}{\sfrq{qa, \z{qa}{b'\ c'}}{L}}
\z{\sfrq{\z{qa}{b''}, \z{qa}{c''}}{K}}{\sfrq{qa, \z{qa}{b''\ c''}}{K}}
$$
$$ \times\ 
\vwp{12}{11}{a; b,c,b'\ , c'\ , b''\ , c''\ , q^{-N}, 
q^{-L}, q^{-K};q\ ,\z{a^{4} q^{N+L+K+4}}{b\ c\ b'\ c'\ b''\ c''}} . 
\eqno(3.5)$$
\subsubsection*{Proof of (3.1)}
For the choice
$$ 
\ez{\aph}{\sfrq{a,q\sqrt{a}, -q\sqrt{a}, d }{r}\ 
q^{r^{2}-r}\ \lrdl{(}{)}{\z{aq^{3}}{d}}^{r}}
{\sfrq{\sqrt{a}, -\sqrt{a}, \z{qa}{d},\ q}{r}}\ ,
\ez{\delta}{\sfrq{b,\ c,\ q^{-M}}{r}}{\sfrq{\z{bcq^{-M}}{a}}{r}}\ ,
$$\\[9pt] 
$$
\ez{\delta'}{\sfrq{b',\ c',\ q^{-N}}{r}}{\sfrq{\z{b'\ c'\ q^{-N}}
{a}}{r}},\  v_{r}= v'_{r}=\z{1}{\sfrq{qa}{r}}, u_{r}=u'_{r}
=\z{q^{r}}{(q;q)_{r}}, 
$$
\\ 
and $t_{r}=w_{r}=1$  in (2.1) and (2.2), the 
\qpss\ and the \sxpf\  can be used to have
$$
\bez{\sfrq{\z{qa}{d}}{n+l}\ q^{n+l}}
{\sfrq{qa}{n+l} \sfrq{\z{qa}{d}, q}{n} \sfrq{\z{qa}{d}, q}{l}}
$$ and 
$$
\gmm_{m}=\z{\sfrq{\z{qa}{b}, \z{qa}{c}}{M}}
{\sfrq{qa, \z{qa}{bc}}{M}}
\z{\sfrq{\z{qa}{b'}, \z{qa}{c'}}{N}}
{\sfrq{qa, \z{qa}{b'\ c'}}{N}}
\z{\sfrq{b, c, b',\ c',\ q^{-M}, q^{-N}}{m}q^{-m^{2}+m}}
{\sfrq{\z{qa}{b}, \z{qa}{c}, \z{qa}{b'}, \z{qa}{c'}, 
aq^{1+M}, aq^{1+N}}{m}}\lrdl{(}{)}{\z{a^{2}\ q^{M+N}}
{b\ c\ b'\ c'}}^{m}
.$$ \\
Substituting these values in (2.3), we obtain the result (3.1).
\subsubsection*{Proof of (3.2)}
For the choice
$$ 
\ez{\aph}{\sfrq{a,q\sqrt{a}, -q\sqrt{a}, \z{qa^{2}}{AB} }{r}\ 
 \lrdl{(}{)}{\z{qAB}{a^{2}}}^{r}}
{\sfrq{\sqrt{a}, -\sqrt{a}, \z{AB}{a},\ q}{r}}\ ,
$$
$$\ez{\delta}{\sfrq{q\sqrt{A}, -q\sqrt{A},\ b,\ c, \z{Aaq^{1+M}}{bc},
\ q^{-M}}{r}}{\sfrq{\sqrt{A}, -\sqrt{A}, \z{qA}{b}, \z{qA}{c}, 
\z{bcq^{-M}}{a}, Aq^{1+M}}{r}}
\ ,
$$\\[9pt] 
$$
\ez{\delta'}{\sfrq{q\sqrt{B}, -q\sqrt{B},\ b',\ c', \z{Baq^{1+N}}{b'\ c'},
\ q^{-N}}{r}}{\sfrq{\sqrt{B}, -\sqrt{B}, \z{qB}{b'}, \z{qB}{c'}, 
\z{b'\ c'\ q^{-N}}{a}, Bq^{1+N}}{r}}
\ ,
$$
$$  v_{r}=\z{\sfrq{A}{r}}{\sfrq{qa}{r}},\  v'_{r}=
\z{\sfrq{B}{r}}{\sfrq{qa}{r}}, 
u_{r}=\z{\sfrq{\z{A}{a}}{r} q^{r}}{(q;q)_{r}},
u'_{r}
=\z{\sfrq{\z{B}{a}}{r} \ q^{r}}{(q;q)_{r}}, 
$$
\\ 
and $t_{r}=w_{r}=1$  in (2.1) and (2.2) and using the 
\qjeps\  we can obtain
$$
\bez{\sfrq{\z{A}{a}}{n-l}\ \sfrq{\z{AB}{a}}{n+l}\ 
\sfrq{A, \z{qa}{B}}{n}\ \sfrq{B, \z{qa}{A}}{l}}
{\sfrq{\z{qa}{B}}{n-l}\ \sfrq{qa}{n+l}\ \sfrq{\z{AB}{a},q}{n}
\ \sfrq{\z{AB}{a},q}{l}}q^{n}\lrdl{(}{)}{\z{AB}{a^{2}}}^{l}
$$ and 
$$\hspace{-3cm}
\gmm_{m}=\z{\sfrq{qA, \z{qa}{b}, \z{qa}{c}, \z{qA}{bc}}{M}}
{\sfrq{qa, \z{qA}{b}, \z{qA}{c}, \z{qa}{bc}}{M}}
\z{\sfrq{qB, \z{qa}{b'}, \z{qa}{c'}, \z{qB}{b'\ c'}}{N}}
{\sfrq{qa, \z{qB}{b'}, \z{qB}{c'}, \z{qa}{b'\ c'}}{N}}
$$
$$\times\ \z{\sfrq{b, c, \z{Aaq^{1+M}}{bc}, q^{-M}, 
b',\ c',\ \z{Baq^{1+N}}{b'\ c'}, q^{-N}}{m}}
{\sfrq{\z{qa}{b}, \z{qa}{c}, \z{bcq^{-M}}{A}, aq^{1+M}, 
\z{qa}{b'}, \z{qa}{c'}, \z{b'\ c'q^{-N}}{B}, aq^{1+N} }{m}}
\lrdl{(}{)}{\z{a^{2}}{AB}}^{m}
.$$ \\
Substituting these values in (2.3), we obtain the result (3.2).
\subsubsection*{Proof of (3.3)}
For the choice
$$ 
\ez{\aph}{\sfrq{a,q\sqrt{a}, -q\sqrt{a}, d }{r}\ q^{r^{2}}\ 
 \lrdl{(}{)}{\z{q^{2}a}{d}}^{r}}
{\sfrq{\sqrt{a}, -\sqrt{a}, \z{qa}{d},\ q}{r}}\ ,
$$
$$\ez{\delta}{\sfrq{q\sqrt{b'}, -q\sqrt{b'}, b, c, 
\z{ab'\ q^{1+N}}{bc}, q^{-N}, \z{b'\ c'}{a}, \z{b'\ q^{-M}}{a}}{r}}
{\sfrq{\sqrt{b'}, -\sqrt{b'}, \z{qb'}{b}, \z{qb'}{c}, 
\z{bcq^{-N}}{a}, b'\ q^{1+N}}{r}}
\ ,
$$\\[9pt] 
$$  v_{r}=\z{1}{\sfrq{qa}{r}}=  v'_{r}, 
u_{r}=\z{ q^{r}}{(q;q)_{r}}= u'_{r},\ 
\ez{w}{\sfrq{b'}{r}}{\sfrq{\z{b'\ c'\ q^{-M}}{a}}{r}}, 
$$
\\ 
$ \delta'_{r}=\sfrq{c',\ q^{-M}}{r}$and $t_{r}=1$ 
 in (2.1) and (2.2), the \sxpf , the \qpss\ and the 
\qjeps\  leads to  
$$
\bez{\sfrq{b',\ \z{qa}{d}}{n+l}}{\sfrq{aq, \z{b'\ c'\ q^{-M}
}{a}}{n+l}} \z{q^{n+l}}{\sfrq{\z{qa}{d},q}{n} \sfrq{\z{qa}{d},q}{l}}
$$ and 
$$\hspace{-3cm}
\gmm_{m}=\z{\sfrq{ \z{qa}{b'}, \z{qa}{c'} }{M}}
{\sfrq{ qa, \z{qa}{b'\ c'}}{M}}
\z{\sfrq{qb', \z{qa}{b}, \z{qa}{c}, \z{qb'}{b\ c}}{N}}
{\sfrq{qa, \z{qb'}{b}, \z{qb'}{c}, \z{qa}{b\ c}}{N}}
$$
$$\times\ \z{\sfrq{b, c, c',\ \z{ab'q^{1+N}}{bc}, q^{-M}, 
 q^{-N}}{m}}
{\sfrq{\z{qa}{b}, \z{qa}{c}, \z{qa}{c'}, \z{bcq^{-N}}{b'},
 aq^{1+M}, aq^{1+N} }{m}}
\lrdl{(}{)}{\z{a}{b'\ c'}}^{m}q^{-m^{2}+mM}
.$$ \\
Substituting these values in (2.3), we obtain the result (3.3).
\subsubsection*{Proof of (3.4)}
For the choice
$$ 
\ez{\aph}{\sfrq{a,q\sqrt{a}, -q\sqrt{a}, \z{a^{2}q^{1+M}}{b'\ c'}
 }{r}\ 
 \lrdl{(}{)}{\z{qb'\ c'}{a^{2}q^{M}}}^{r}}
{\sfrq{\sqrt{a}, -\sqrt{a}, \z{b'\ c'\ q^{-M}}{a},\ q}{r}}\ ,
$$
$$\ez{\delta}{\sfrq{q\sqrt{b'}, -q\sqrt{b'}, b, c, 
\z{ab'\ q^{1+N}}{bc}, q^{-N}, \z{b'\ c'}{a}, 
\z{b'\ q^{-M}}{a}}{r}}
{\sfrq{\sqrt{b'}, -\sqrt{b'}, \z{qb'}{b}, \z{qb'}{c}, 
\z{bcq^{-N}}{a}, b'\ q^{1+N}}{r}}
\ ,
$$\\[9pt]
$$
\delta'_{r}=\z{{\sfrq{q\sqrt{A}, -q\sqrt{A},  c',\  
q^{-M}}{r}}}
{{\sfrq{\sqrt{A}, -\sqrt{A},  \z{qA}{c'}, 
Aq^{1+M}}{r}}}\ , w_{r}=\z{\sfrq{b'}{r}}
{\sfrq{\z{b'\ c'\ q^{-M}}{a}}{r}},
\ez{t}{\sfrq{\z{b'}{A}}{r}}{\sfrq{\z{b'\ c'\ q^{-M}}{Aa}}{r}},$$ 
\\[9pt] $$ v_{r}=\z{\sfrq{\z{b'\ c'\ q^{-M}}{A}}{r}}{\sfrq{qa}{r}},\  
v'_{r}= \z{\sfrq{A}{r}}{\sfrq{qa}{r}}, 
u_{r}
=\z{\sfrq{\z{b'\ c'\ q^{-M}}{Aa}}{r} \ q^{r}}{(q;q)_{r}} 
\mbox{ and }u'_{r}=\z{\sfrq{\z{A}{a}}{r} q^{r}}{(q;q)_{r}},
$$
 \\ 
  in (2.1) and (2.2),  the \qjeps\  gives rise to  \\
$$
\bez{\sfrq{\z{a}{A}}{n-l}\ \sfrq{b'}{n+l}\ 
\sfrq{\z{b'\ c'\ q^{-M}}{A}, \z{qa}{A}}{n}\ 
\sfrq{A, \z{aAq^{1+M}}{b'\ c'}}{l}}
{\sfrq{\z{qa}{A}}{n-l}\ \sfrq{qa}{n+l}\ 
\sfrq{\z{b'\ c'\ q^{-M}}{A},q}{n}
\ \sfrq{\z{b'\ c'\ q^{-M}}{A},q}{l}}q^{n}
\lrdl{(}{)}{\z{b'\ c'\ q^{-M}}{a^{2}}}^{l}
$$ and 
$$\hspace{-3cm}
\gmm_{m}=\z{\sfrq{qA, \z{qa}{b'}, \z{qa}{c'}, \z{qA}{b'\ c'}}{M}}
{\sfrq{qa, \z{qA}{b'}, \z{qA}{c'}, \z{qa}{b'\ c'}}{M}}
\z{\sfrq{qb', \z{qa}{b}, \z{qa}{c}, \z{qb'}{b c}}{N}}
{\sfrq{qa, \z{qb'}{b}, \z{qb'}{c}, \z{qa}{b c}}{N}}
$$
$$\times\ \z{\sfrq{b, c, c',\ \z{b'aq^{1+N}}{bc}, q^{-M}, 
 q^{-N}}{m}}
{\sfrq{\z{qa}{b}, \z{qa}{c}, \z{qa}{c'}, \z{bcq^{-N}}{b'},
 aq^{1+M},  aq^{1+N} }{m}}
\lrdl{(}{)}{\z{a}{b'}}^{m}
.$$ \\[9pt]
Substituting these values in (2.3), we obtain the result (3.4).
\subsubsection*{Proof of (3.5)}
For the choice
$$ 
\ez{\aph}{\sfrq{a,q\sqrt{a}, -q\sqrt{a} }{r}\ 
a^{r}q^{4r}\ q^{3\binom{r}{2}}}
{\sfrq{\sqrt{a}, -\sqrt{a}, \ q}{r}}\ ,
\ez{\delta}{\sfrq{b,\ c,\ q^{-N}}{r}}{\sfrq{\z{bcq^{-N}}{a}}{r}}\ ,
$$\\[9pt] 
$$
\ez{\delta'}{\sfrq{b',\ c',\ q^{-L}}{r}}{\sfrq{\z{b'\ c'\ q^{-L}}
{a}}{r}},\  
\ez{\delta''}{\sfrq{b'',\ c'',\ q^{-K}}{r}}{\sfrq
{\z{b''\ c''\ q^{-K}}{a}}{r}},
$$
$$ v_{r}= v'_{r}= v''_{r}=\z{1}{\sfrq{qa}{r}}, \mbox{  and  } 
\ u_{r}=u'_{r} = u''_{r}= \z{q^{r}}{(q;q)_{r}}, 
$$
\\ 
  in (2.6) and (2.7), the 
\qpss\ and the \sxpf\  readily yields 
$$
\bezk{\sfrq{qa}{n+l+k}\ q^{n+l+k}}
{\sfrq{qa}{n+l}\sfrq{qa}{n+k}\sfrq{qa}{l+k} 
\sfrq{q}{n} \sfrq{q}{l}\sfrq{q}{k}}
$$ and 
$$
\gmm_{m}=\z{\sfrq{\z{qa}{b}, \z{qa}{c}}{N}}
{\sfrq{qa, \z{qa}{bc}}{N}}
\z{\sfrq{\z{qa}{b'}, \z{qa}{c'}}{L}}
{\sfrq{qa, \z{qa}{b'\ c'}}{L}}
\z{\sfrq{\z{qa}{b''}, \z{qa}{c''}}{K}}
{\sfrq{qa, \z{qa}{b''\ c''}}{K}}$$
$$ \z{\sfrq{b, c, b',\ c',\ b'',\ c'',\ q^{-N}, q^{-L}, q^{-K}}
{m}(-1)^{m} q^{-3\binom{m}{2}}}
{\sfrq{\z{qa}{b}, \z{qa}{c}, \z{qa}{b'}, \z{qa}{c'}, 
\z{qa}{b''}, \z{qa}{c''},aq^{1+N}, aq^{1+L}, aq^{1+K}}{m}}
\lrdl{(}{)}{\z{a^{3}\ q^{N+L+K}} {b\ c\ b'\ c'\ b''\ c''}}^{m}
.$$ \\[9pt]
Substituting these values in (2.8), we obtain the result (3.5).
\section{First Bailey type lemma and applications to \rr\ type \id\ }
\subsection{First Bailey Type Lemma or FBTL}
\subsubsection*{Theorem 4.1}
If for $ n, l \geq 0 $ 
$$
\btt_{(n, l)} = \asmm{m=0}{min. (n, l)}
\z{\aph_{m}}{(q; q)_{n-m}(q; q)_{l-m}(aq;q)_{n+m}(aq;q)_{l+m}}
 , \eqno(4.1)$$ then 
 $$
 \btt'_{(n, l)} = \asmm{m=0}{min. (n, l)}
\z{\aph'_{m}}{(q; q)_{n-m}(q; q)_{l-m}(aq;q)_{n+m}(aq;q)_{l+m}}
 , \eqno(4.2) 
 $$ where 
 $$
 \aph'_{m}= \z{\sfrq{b, c, b',\ c'}{m} \lrdl{(}{)}{\z{a^{2}q^{2}}
{bcb'\ c'}}^{m} \ \aph_{m}}
{\sfrq{\z{qa}{b}, \z{qa}{c}, \z{qa}{b'}, \z{qa}{c'}}{m}}
 \eqno(4.3)$$ and 
 $$
 \btt'_{(M, N)}= \asmm{n, l =0}{\infty} 
 \z{\sfrq{b,c}{n} \sfrq{\z{qa}{bc}}{M-n}\lrdl{(}{)}
{\z{qa}{bc}}^{n}}
{\sfrq{\z{qa}{b}, \z{qa}{c}}{M} (q; q)_{M-n}}
\z{\sfrq{b',c'}{l} \sfrq{\z{qa}{b'\ 'c}}{N-l}\lrdl{(}{)}
{\z{qa}{b'\ c'}}^{l}\ \btt_{(n, l)}}
{\sfrq{\z{qa}{b}, \z{qa}{c}}{N} (q; q)_{N-l}} . 
 \eqno(4.4)$$
\textsc{Remark:} A pair of sequences $( \aph_{m}, \btt_{(n, l)})$ 
related by (4.1) may be called as a ``first Bailey type pair" or 
``FBTP" and the Theorem 4.1 may be rephrased as: If 
$( \aph_{m}, \btt_{(n, l)})$ is a FBTP, so is 
 $( \aph'_{m}, \btt'_{(n, l)})$ where this new FBTP 
 is given by (4.3) and (4.4).
\subsubsection*{Proof: } To prove FBTL, we apply 
Bailey type transform (2.3) with the choice:
\begin{equation*}
\ez{\delta}{\sfrq{b,\ c,\ q^{-M}}{r}q^{r}}{\sfrq{\z{bcq^{-M}}{a}}{r}}\ ,
\ez{\delta'}{\sfrq{b',\ c',\ q^{-N}}{r}q^{r}}{\sfrq{\z{b'\ c'\ q^{-N}}
{a}}{r}},\end{equation*}
\begin{equation*} 
v_{r}= v'_{r}=\z{1}{\sfrq{qa}{r}}, u_{r}=u'_{r}
=\z{1}{(q;q)_{r}},
\end{equation*}
\\ 
and $t_{r}=w_{r}=1$ then as in the proof of (3.1), we get
\begin{multline*} \hspace{-1cm}
\gmm_{m}=\z{\sfrq{\z{qa}{b}, \z{qa}{c}}{M}}
{\sfrq{qa, \z{qa}{bc}}{M}}
\z{\sfrq{\z{qa}{b'}, \z{qa}{c'}}{N}}
{\sfrq{qa, \z{qa}{b'\ c'}}{N}}\\[9pt]
\times\ \z{\sfrq{b, c, b',\ c',\ q^{-M}, q^{-N}}{m}q^{-m^{2}+m}}
{\sfrq{\z{qa}{b}, \z{qa}{c}, \z{qa}{b'}, \z{qa}{c'}, 
aq^{1+M}, aq^{1+N}}{m}}\lrdl{(}{)}{\z{a^{2}\ q^{2+M+N}}
{b\ c\ b'\ c'}}^{m}
.
\end{multline*} 
\\
Now we can prove (4.2) as shown below:
\begin{multline*} \hspace{-1cm}
\asmm{m=0}{min. (M, N)}
\z{\aph'_{m}}{(q; q)_{M-m}(q; q)_{N-m}(aq;q)_{M+m}(aq;q)_{N+m}}\\
 \\ =
  \asmm{m=0}{min. (M, N)}
\z{\sfrq{b, c, b',\ c'}{m} \lrdl{(}{)}{\z{a^{2}q^{2}}
{bcb'\ c'}}^{m} \ \aph_{m}}
{\sfrq{\z{qa}{b}, \z{qa}{c}, \z{qa}{b'}, \z{qa}{c'}}{m}
(q; q)_{M-m}(q; q)_{N-m}(aq;q)_{M+m}(aq;q)_{N+m}} \\ \\ 
=\z{1}
{\sfrq{qa, q}{M}}
\z{1}
{\sfrq{qa, q}{N}} \\ \\ 
\times\ \asmm{m=0}{min. (M, N)}
\z{\sfrq{b, c, b',\ c',\ q^{-M}, q^{-N}}{m}q^{-m^{2}+m}}
{\sfrq{\z{qa}{b}, \z{qa}{c}, \z{qa}{b'}, \z{qa}{c'}, 
aq^{1+M}, aq^{1+N}}{m}}\lrdl{(}{)}{\z{a^{2}\ q^{M+N}}
{b\ c\ b'\ c'}}^{m} \aph_{m} \\ \\ 
=\z{\sfrq{ \z{qa}{bc}}{M}}
{\sfrq{\z{qa}{b}, \z{qa}{c}, q}{M}}
\z{\sfrq{ \z{qa}{b'\ c'}}{N}}
{\sfrq{\z{qa}{b'}, \z{qa}{c'}, q}{N}}
\asmm{m=0}{min. (M, N)} \gmm_{m}\ \aph_{m} \\ \\ 
=\z{\sfrq{ \z{qa}{bc}}{M}}
{\sfrq{\z{qa}{b}, \z{qa}{c}, q}{M}}
\z{\sfrq{ \z{qa}{b'\ c'}}{N}}
{\sfrq{\z{qa}{b'}, \z{qa}{c'}, q}{N}}
\asmm{n=0}{M}\asmm{l=0}{N} \btt_{(n, l)}\ \delta_{n}
\delta'_{l} \hspace{1.5cm} \mbox{(by Theorem 4.1)} 
\end{multline*}
\\
\begin{multline*}
=\z{\sfrq{ \z{qa}{bc}}{M}}
{\sfrq{\z{qa}{b}, \z{qa}{c}, q}{M}}
\z{\sfrq{ \z{qa}{b'\ c'}}{N}}
{\sfrq{\z{qa}{b'}, \z{qa}{c'}, q}{N}}\\ \\
\times\ \asmm{n=0}{M}\asmm{l=0}{N} \z{\sfrq{b,\ c,\ q^{-M}}{n}q^{n} 
\ \sfrq{b',\ c',\ q^{-N}}{l}q^{l}\ \btt_{(n, l)}}
{\sfrq{\z{bcq^{-M}}{a}}{n}\ \sfrq{\z{b'\ c'\ q^{-N}}{a}}{l}}
 \\ \\ 
= \btt'_{(M, N)}. \\
\end{multline*}
The last line follows by manipulating \qsf , which reduces the 
previous expression to (4.4). $\quad \quad \square $
\\
  \\
 Now to derive \rr\ type \id , we need  
the following result obtained by substituting the values of 
$ \aph'_{m}$ and $ \btt'_{(n, l)}$ from (4.3) and (4.4) into (4.2),
 viz. 
\begin{multline} \hspace{-1.1cm}
\asmm{n, l =0}{\infty} 
 \z{\sfrq{b,c}{n} \sfrq{\z{qa}{bc}}{M-n}\lrdl{(}{)}
{\z{qa}{bc}}^{n}}
{\sfrq{\z{qa}{b}, \z{qa}{c}}{M} (q; q)_{M-n}}
\z{\sfrq{b',c'}{l} \sfrq{\z{qa}{b'\ 'c}}{N-l}\lrdl{(}{)}
{\z{qa}{b'\ c'}}^{l}\ \btt_{(n, l)}}
{\sfrq{\z{qa}{b}, \z{qa}{c}}{N} (q; q)_{N-l}} \\ \\
\hspace{-1cm} =
  \asmm{m=0}{min. (M, N)}
\z{\sfrq{b, c, b',\ c'}{m} \lrdl{(}{)}{\z{a^{2}q^{2}}
{bcb'\ c'}}^{m} \ \aph_{m}}
{\sfrq{\z{qa}{b}, \z{qa}{c}, \z{qa}{b'}, \z{qa}{c'}}{m}
(q; q)_{M-m}(q; q)_{N-m}(aq;q)_{M+m}(aq;q)_{N+m}}
\tag{4.5}
\end{multline} 
Now taking $ b,c, b',\ c',\ M, N \rightarrow \infty $ in (4.5), 
we get 
$$
  \asmm{n, l =0}{\infty} a^{n+l}\ 
q^{n^{2}+l^{2}}\ \btt_{(n, l)} = \z{1}{\sfrq{aq}{\infty}^{2}} \ \ 
\asmm{m=0}{\infty} a^{2m}\ q^{2m^{2}}\ \aph_{m} ,
 \eqno(4.6)$$
for any FBTP $( \aph_{m}, \btt_{(n, l)})$. 
Now, we shall make use of the following two ``FBTPs" deduced by us : 
\begin{align*} \begin{split} 
&\aph_{m}= \z{\sfrq{a}{m} (1-aq^{2m})}{\sfrq{q}{m} (1-a)} 
(-1)^{m} a^{m}\ q^{\frac{1}{2} (3m^{2}-m)} ,  \\[9pt] 
&\btt_{(n, l)}= \z{1}{\sfrq{qa}{n+l} \sfrq{q}{n} \sfrq{q}{l}} , 
\end{split} \tag{4.7} 
\end{align*} and 
\begin{align*} \begin{split} 
&\aph_{m}= \z{\sfrq{a}{m} (1-aq^{2m})}{\sfrq{q}{m} (1-a)} 
(-1)^{m} \ q^{\frac{1}{2} (m^{2}-m)} ,  \\[9pt] 
&\btt_{(n, l)}= \z{q^{n l}}{\sfrq{qa}{n+l} \sfrq{q}{n} \sfrq{q}{l}} . 
\end{split} \tag{4.8} 
\end{align*}
The fact that ``FBTPs" $( \aph_{m}, \btt_{(n, l)})$ given by (4.7) and (4.8) 
satisfy (4.1) may be verified by substituting them in (4.1) 
and then appealing to \sxpf\ and taking $d\rightarrow \infty$ or 
$d\rightarrow 0$. \\[9pt]   
\subsection{New double series \rr\ type \id\ and 
corresponding infinite families  :} To derive \rr\ type \id , 
we insert the  ``FBTPs" (4.7) and (4.8) in (4.6) to get (4.9) and (4.10), 
respectively, as given below: 
\begin{multline}
\asmm{n, l =0}{\infty} \z{ a^{n+l}\ 
q^{n^{2}+l^{2}}}{\sfrq{qa}{n+l} \sfrq{q}{n} \sfrq{q}{l}}\\[9pt]
= \z{1}{\sfrq{aq}{\infty}^{2}}
\asmm{m=0}{\infty} \z{\sfrq{a}{m} (1-aq^{2m})}{\sfrq{q}{m} (1-a)} 
 (-1)^{m} a^{3m}\ q^{\frac{1}{2} (7m^{2}-m)} ,  
\tag{4.9}
\end{multline}
\begin{multline}
\asmm{n, l =0}{\infty} \z{ a^{n+l}\ 
q^{n^{2}+l^{2}+nl}}{\sfrq{qa}{n+l} \sfrq{q}{n} \sfrq{q}{l}}\\[9pt]
= \z{1}{\sfrq{aq}{\infty}^{2}}
\asmm{m=0}{\infty} \z{\sfrq{a}{m} (1-aq^{2m})}{\sfrq{q}{m} (1-a)} 
 (-1)^{m} a^{2m}\ q^{\frac{1}{2} (5m^{2}-m)}  , 
\tag{4.10}
\end{multline}
and then setting $ a=1 $ and $ a=q $ in (4.9) and (4.10) and
 using Jacobi triple product identity, we obtain following four 
 new double series \rr\ type \id :
 $$
\asmm{n, l =0}{\infty} \z{  
q^{n^{2}+l^{2}}}{\sfrq{q}{n+l} \sfrq{q}{n} \sfrq{q}{l}}
= \z{1}{\sfrq{q}{\infty}} \spdt{0, \pm 3 \text{ (mod\ 7) }} , 
 \eqno(4.11)
 $$
 $$
\asmm{n, l =0}{\infty} \z{  
q^{n^{2}+l^{2}+n+l}}{\sfrq{q^{2}}{n+l} \sfrq{q}{n} \sfrq{q}{l}}
= \z{1}{\sfrq{q^{2}}{\infty}} \spdt{0, \pm 1 \text{(mod\ 7)}} , 
 \eqno(4.12)
 $$ \\[9pt] 
$$
\asmm{n, l =0}{\infty} \z{  
q^{n^{2}+l^{2}+nl}}{\sfrq{q}{n+l} \sfrq{q}{n} \sfrq{q}{l}}
= \z{1}{\sfrq{q}{\infty}} \spdt{0, \pm 2 \text{(mod\ 5)}} , 
 \eqno(4.13)
 $$ \\ and 
 $$
\asmm{n, l =0}{\infty} \z{  
q^{n^{2}+l^{2}+n+l+nl}}{\sfrq{q^{2}}{n+l} \sfrq{q}{n} \sfrq{q}{l}}
= \z{1}{\sfrq{q^{2}}{\infty}} 
\spdt{0, \pm 1 \text{(mod\ 5)}} . 
 \eqno(4.14)
 $$
\\ 
Further, like classical Bailey lemma, the idea of repeated 
application of FBTL may be expressed in the concept
 of the ``first Bailey type chain." Given a FBTP 
$ (\aph_{m}, \btt_{(n, l)}) $ we can by FBTL produce new FBTP 
$ (\aph'_{m}, \btt'_{(n, l)}) $ defined by (4.3) and (4.4). 
From $(\aph'_{m}, \btt'_{(n, l)}) $ we can create 
$ (\aph''_{m}, \btt''_{(n, l)})$ merely by applying 
FBTL with $ (\aph'_{m}, \btt'_{(n, l)}) $ as 
initial FBTP. In this way we create a 
sequence of FBTPs:
$$
(\aph_{m}, \btt_{(n, l)})\rightarrow (\aph'_{m}, \btt'_{(n, l)})\rightarrow 
(\aph''_{m}, \btt''_{(n, l)})\rightarrow (\aph'''_{m}, \btt'''_{(n, l)})\rightarrow 
\cdots ;
$$
and call it ``first Bailey type chain." Using this chain 
we can establish:
\subsubsection*{Theorem 4.2}
\begin{align*} \hspace{-1cm}
\asmm{m \geq 0}{} & \z{\sfrq{b_{1}}{m}\sfrq{c_{1}}{m} 
\sfrq{b_{2}}{m}\sfrq{c_{2}}{m} \cdots 
\sfrq{b_{k}}{m}\sfrq{c_{k}}{m}}
{\sfrq{\z{aq}{b_{1}}}{m}\sfrq{\z{aq}{c_{1}}}{m} 
\sfrq{\z{aq}{b_{2}}}{m}\sfrq{\z{aq}{c_{2}}}{m} \cdots 
\sfrq{\z{aq}{b_{k}}}{m}\sfrq{\z{aq}{c_{k}}}{m}} 
\\ \\
& \times\ \z{\sfrq{b'_{1}}{m}\sfrq{c'_{1}}{m} 
\sfrq{b'_{2}}{m}\sfrq{c'_{2}}{m} \cdots 
\sfrq{b'_{k}}{m}\sfrq{c'_{k}}{m}}
{\sfrq{\z{aq}{b'_{1}}}{m}\sfrq{\z{aq}{c'_{1}}}{m} 
\sfrq{\z{aq}{b'_{2}}}{m}\sfrq{\z{aq}{c'_{2}}}{m} \cdots 
\sfrq{\z{aq}{b'_{k}}}{m}\sfrq{\z{aq}{c'_{k}}}{m}} \\ \\ 
& \times\ \z{\sfrq{q^{-M}}{m}\sfrq{q^{-N}}{m}}
{\sfrq{aq^{1+M}}{m}\sfrq{aq^{1+N}}{m}} 
\lrdl{(}{)}{\frac{a^{2k}\ q^{2k+M+N}}
{b_{1}c_{1} \cdots b_{k}c_{k} b'_{1}c'_{1} 
\cdots b'_{k}c'_{k}}}^{m}
 q^{-m^{2}+m}\ \aph_{m} \\ \\ 
& =\z{\sfrq{aq}{M} \sfrq{\z{aq}{b_{k}c_{k}}}{M}}
{\sfrq{\z{aq}{b_{k}}}{M}\sfrq{\z{aq}{c_{k}}}{M}} 
\z{\sfrq{aq}{N} \sfrq{\z{aq}{b'_{k}c'_{k}}}{N}}
{\sfrq{\z{aq}{b'_{k}}}{N}\sfrq{\z{aq}{c'_{k}}}{N}}
\\ \\ & \times\ 
\asmm{\substack{n_{k} \geq \cdots \geq n_{1} \geq 0 
\\ \\  l_{k} \geq \cdots \geq l_{1} \geq 0}}{} 
\z{\sfrq{b_{k}}{n_{k}}\sfrq{c_{k}}{n_{k}} 
 \cdots 
\sfrq{b_{1}}{n_{1}}\sfrq{c_{1}}{n_{1}}} 
{\sfrq{q}{n_{k}-n_{k-1}}\sfrq{q}{n_{k-1}-n_{k-2}} 
\cdots \sfrq{q}{n_{2}-n_{1}}}\\ \\ & \times\
\z{\sfrq{b'_{k}}{l_{k}}\sfrq{c'_{k}}{l_{k}} 
 \cdots 
\sfrq{b'_{1}}{l_{1}}\sfrq{c'_{1}}{l_{1}}}
{\sfrq{q}{l_{k}-l_{k-1}}\sfrq{q}{l_{k-1}-l_{k-2}} 
\cdots \sfrq{q}{l_{2}-l_{1}}}
\z{\sfrq{q^{-M}}{n_{k}}}
{\sfrq{\z{b_{k}\ c_{k}q^{-M}}{a}}{n_{k}}}
\z{\sfrq{q^{-N}}{l_{k}}}
{\sfrq{\z{b'_{k}\ c'_{k}q^{-N}}{a}}{l_{k}}}
\\ \\ & \times\ 
\z{ \sfrq{\z{aq}{b_{k-1}\ c_{k-1}}}{n_{k}-n_{k-1}}
\cdots \sfrq{\z{aq}{b_{1}\ c_{1}}}{n_{2}-n_{1}}
}
{\sfrq{\z{aq}{b_{k-1}}}{n_{k}}\sfrq{\z{aq}{c_{k-1}}}{n_{k}} 
\cdots \sfrq{\z{aq}{b_{1}}}{n_{2}}\sfrq{\z{aq}{c_{1}}}{n_{2}}
} \end{align*} 
\begin{align*} & \times\ 
\z{ \sfrq{\z{aq}{b'_{k-1}\ c'_{k-1}}}{l_{k}-l_{k-1}}
\cdots \sfrq{\z{aq}{b'_{1}\ c'_{1}}}{l_{2}-l_{1}}
}
{\sfrq{\z{aq}{b'_{k-1}}}{l_{k}}\sfrq{\z{aq}{c'_{k-1}}}{l_{k}} 
\cdots \sfrq{\z{aq}{b'_{1}}}{l_{2}}\sfrq{\z{aq}{c'_{1}}}{l_{2}}
} \\ \\ & \times\ 
q^{n_{1}+n_{2}+ \cdots +n_{k}+l_{1}+l_{2}+ \cdots +l_{k}} \ 
a^{n_{1}+n_{2}+ \cdots +n_{k-1}+l_{1}+l_{2}+ \cdots +l_{k-1}} 
\\ \\ & \times\ 
\lrdl{(}{)}{b_{k-1}c_{k-1}}^{-n_{k-1}} \cdots 
\lrdl{(}{)}{b_{1}c_{1}}^{-n_{1}} 
\lrdl{(}{)}{b'_{k-1}c'_{k-1}}^{-l_{k-1}} \cdots 
\lrdl{(}{)}{b'_{1}c'_{1}}^{-l_{1}}\ \btt_{(n_{1}, l_{1})} . 
\tag{4.15} 
\end{align*}
\subsubsection*{Proof :}
Equation (4.15) is the result of  a $k$-fold iteration of FBTL 
with the choice $ b=b_{i}$, $c=c_{i}, b'=b'_{i}$, $c'=c'_{i} $ at the $i$ th step. 
$ \square $ 
\subsubsection*{Theorem 4.3}
If $ (\aph_{m}, \btt_{(n, l)}) $ is a FBTP (i.e. related by (4.1)),
 then 
\begin{multline*} \hspace{-1cm} 
\asmm{\substack{n_{k} \geq \cdots \geq n_{1} \geq 0 
\\ \\  l_{k} \geq \cdots \geq l_{1} \geq 0}}{}  
\z{ a^{n_{1}+n_{2}+ \cdots +n_{k}+l_{1}+l_{2}+ \cdots +l_{k}} 
q^{n^{2}_{1}+n^{2}_{2}+ \cdots +n^{2}_{k}
+l^{2}_{1}+l^{2}_{2}+ \cdots +l^{2}_{k}}\ \btt_{(n_{1}, l_{1})}
}
{\sfrq{q}{n_{k}-n_{k-1}} 
\cdots \sfrq{q}{n_{2}-n_{1}}
\sfrq{q}{l_{k}-l_{k-1}} 
\cdots \sfrq{q}{l_{2}-l_{1}}}\\[9pt] 
 = \z{1}{\sfrq{aq}{\infty}^{2}}
\asmm{m=0}{\infty} q^{2km^{2}}\ a^{2km}\ \aph_{m} . 
\tag{4.16}
\end{multline*} 
\subsubsection*{Proof :}
Let $M, N, b_{1},  \cdots, b_{k}, c_{1},  \cdots c_{k},
b'_{1},  \cdots, b'_{k}, c'_{1},  \cdots c'_{k} $ 
all tend to infinity in Theorem 4.2. $ \square$  \\[9pt]  
Theorem 4.3 may now be used to embed the double series \rr\ 
type \id\ (4.11) - (4.14) in an infinite family of such identities : 
\begin{multline*} \hspace{-1.3cm} 
\asmm{\substack{n_{k} \geq \cdots \geq n_{1} \geq 0 
\\ \\  l_{k} \geq \cdots \geq l_{1} \geq 0}}{}  
\z{ 
q^{n^{2}_{1}+n^{2}_{2}+ \cdots +n^{2}_{k}
+l^{2}_{1}+l^{2}_{2}+ \cdots +l^{2}_{k}}
}
{\sfrq{q}{n_{k}-n_{k-1}} 
\cdots \sfrq{q}{n_{2}-n_{1}}
\sfrq{q}{l_{k}-l_{k-1}} 
\cdots \sfrq{q}{l_{2}-l_{1}} }\\ \\
\times\ \z{1}
{\sfrq{q}{n_{1}+l_{1}}\sfrq{q}{n_{1}}\sfrq{q}{l_{1}}} 
= \z{1}{\sfrq{q}{\infty}}
 \spdt{0, \pm (2k+1) \text{(mod\ 4k+3)}} , 
\tag{4.17}
\end{multline*} \\  
\begin{multline*} \hspace{-1.3cm} 
\asmm{\substack{n_{k} \geq \cdots \geq n_{1} \geq 0 
\\ \\  l_{k} \geq \cdots \geq l_{1} \geq 0}}{}  
\z{ 
q^{n^{2}_{1}+n^{2}_{2}+ \cdots +n^{2}_{k}
+l^{2}_{1}+l^{2}_{2}+ \cdots +l^{2}_{k}+
n_{1}+n_{2}+ \cdots +n_{k}
+l_{1}+l_{2}+ \cdots +l_{k}}
}
{\sfrq{q}{n_{k}-n_{k-1}} 
\cdots \sfrq{q}{n_{2}-n_{1}}
\sfrq{q}{l_{k}-l_{k-1}} 
\cdots \sfrq{q}{l_{2}-l_{1}} }\\ \\ 
 \times\ \z{1}
{\sfrq{q^{2}}{n_{1}+l_{1}}
\sfrq{q}{n_{1}}\sfrq{q}{l_{1}}}
= \z{1}{\sfrq{q^{2}}{\infty}}
 \spdt{0, \pm 1 \text{(mod\ 4k+3)}} , 
\tag{4.18}
\end{multline*} \\ 
\begin{multline*} \hspace{-1.3cm} 
\asmm{\substack{n_{k} \geq \cdots \geq n_{1} \geq 0 
\\ \\  l_{k} \geq \cdots \geq l_{1} \geq 0}}{}  
\z{ 
q^{n^{2}_{1}+n^{2}_{2}+ \cdots +n^{2}_{k}
+l^{2}_{1}+l^{2}_{2}+ \cdots +l^{2}_{k}+n_{1}l_{1}}
}
{\sfrq{q}{n_{k}-n_{k-1}} 
\cdots \sfrq{q}{n_{2}-n_{1}}
\sfrq{q}{l_{k}-l_{k-1}} 
\cdots \sfrq{q}{l_{2}-l_{1}} }\\ \\ 
\z{1}{\sfrq{q}{n_{1}+l_{1}}
\sfrq{q}{n_{1}}\sfrq{q}{l_{1}}} = \z{1}{\sfrq{q}{\infty}}
 \spdt{0, \pm 2k \text{(mod\ 4k+1)}} 
\tag{4.19}
\end{multline*} \\ and 
\begin{multline*} \hspace{-1.3cm} 
\asmm{\substack{n_{k} \geq \cdots \geq n_{1} \geq 0 
\\ \\  l_{k} \geq \cdots \geq l_{1} \geq 0}}{}  
\z{ 
q^{n^{2}_{1}+n^{2}_{2}+ \cdots +n^{2}_{k}
+l^{2}_{1}+l^{2}_{2}+ \cdots +l^{2}_{k}+
n_{1}+n_{2}+ \cdots +n_{k}
+l_{1}+l_{2}+ \cdots +l_{k}+n_{1}l_{1}}
}
{\sfrq{q}{n_{k}-n_{k-1}} 
\cdots \sfrq{q}{n_{2}-n_{1}}
\sfrq{q}{l_{k}-l_{k-1}} 
\cdots \sfrq{q}{l_{2}-l_{1}} }\\ \\ 
\z{1}{\sfrq{q^{2}}{n_{1}+l_{1}}
\sfrq{q}{n_{1}}\sfrq{q}{l_{1}}}= \z{1}{\sfrq{q^{2}}{\infty}}
\spdt{0, \pm 1 \text{(mod\ 4k+1)}} . 
\tag{4.20}
\end{multline*}
Equation (4.17) follows by using FBTP defined by (4.7) in (4.16) 
and then invoking to Jacobi triple product identity after 
taking $a=1$. For (4.18), we follow the same procedure with $a=q$.
In the same way we establish (4.19) and (4.20) by using FBTP defined by (4.8).
\section{Second Bailey type lemma and applications to \rr\ type \id\ }
\subsection{Second Bailey Type Lemma or SBTL}
\subsubsection*{Theorem 5.1}
If for $ n, l, k \geq 0 $ 
\begin{equation*}\hspace{-1cm}
\btt_{(n, l ,k)}\end{equation*} 
$$ = \asmm{m=0}{min. (n, l, k)}
\z{\aph_{m}}{(q; q)_{n-m}(q; q)_{l-m}(q; q)_{k-m}(aq;q)_{n+m}
(aq;q)_{l+m}(aq;q)_{k+m}}
 , \eqno(5.1)$$ then 
 \begin{equation*} \hspace{-1cm}
\btt'_{(n, l ,k)}\end{equation*} 
$$
 = \asmm{m=0}{min. (n, l, k)}
\z{\aph'_{m}}{(q; q)_{n-m}(q; q)_{l-m}(q;q)_{k-m}(aq;q)_{n+m}
(aq;q)_{l+m}(aq;q)_{k+m}}
 , \eqno(5.2) 
 $$ where 
 $$
 \aph'_{m}= \z{\sfrq{b, c, b',\ c',\ b'',\ c'',}{m} 
\lrdl{(}{)}{\z{a^{3}q^{3}}
{b\ c\ b'\ c'\ b''\ c''}}^{m} \ \aph_{m}}
{\sfrq{\z{qa}{b}, \z{qa}{c}, \z{qa}{b'}, \z{qa}{c'}
\z{qa}{b''}, \z{qa}{c''}}{m}}
 \eqno(5.3)$$ and 
 \begin{multline} \hspace{-1cm}
\btt'_{(N, L, K)}= \asmm{n, l, k =0}{\infty} 
 \z{\sfrq{b,c}{n} \sfrq{\z{qa}{bc}}{N-n}\lrdl{(}{)}
{\z{qa}{bc}}^{n}}
{\sfrq{\z{qa}{b}, \z{qa}{c}}{N} (q; q)_{N-n}}
\z{\sfrq{b',c'}{l} \sfrq{\z{qa}{b'\ c'}}{L-l}\lrdl{(}{)}
{\z{qa}{b'\ c'}}^{l} }
{\sfrq{\z{qa}{b'}, \z{qa}{c'}}{L} (q; q)_{L-l}} \\ \\
\times\ \z{\sfrq{b'',c''}{k} \sfrq{\z{qa}{b''\ c''}}{K-k}\lrdl{(}{)}
{\z{qa}{b''\ c''}}^{k}\ \btt_{(n, l, k)}}
{\sfrq{\z{qa}{b''}, \z{qa}{c''}}{K} (q; q)_{K-k}} . \tag{5.4}
\end{multline}
\textsc{Remark:} A pair of sequences $( \aph_{m}, \btt_{(n, l, k)})$ 
related by (5.1) may be called as a ``second Bailey type pair" or 
``SBTP" and the Theorem 5.1 may be rephrased as: If 
$( \aph_{m}, \btt_{(n, l, k)})$ is a SBTP, so is 
 $( \aph'_{m}, \btt'_{(n, l, k)})$ where this new SBTP 
 is given by (5.3) and (5.4).
\subsubsection*{Proof: } To prove SBTL, we apply 
Bailey type transform (2.8) with the choice:
\begin{equation*}
\ez{\delta}{\sfrq{b,\ c,\ q^{-N}}{r}q^{r}}{\sfrq{\z{bcq^{-N}}{a}}{r}}\ ,
\ez{\delta'}{\sfrq{b',\ c',\ q^{-L}}{r}q^{r}}{\sfrq{\z{b'\ c'\ q^{-L}}
{a}}{r}}, 
\ez{\delta''}{\sfrq{b'',\ c'',\ q^{-K}}{r}q^{r}}{\sfrq{\z{b''\ c''\ q^{-K}}
{a}}{r}}
,\end{equation*} and
\begin{equation*} 
v_{r}= v'_{r}= v''_{r}=\z{1}{\sfrq{qa}{r}}, u_{r}=u'_{r}
=u''_{r}=\z{1}{(q;q)_{r}},
\end{equation*}
\\ 
then following the proof of Theorem 4.1 we can also prove Theorem 5.1 .
\\[9pt]  
Now to derive \rr\ type \id , we need  
the following result obtained by substituting the values of 
$ \aph'_{m}$ and $ \btt'_{(n, l, k)}$ from (5.3) and (5.4) into (5.2),
 {\it viz.,} 
\begin{multline} \hspace{-1.1cm}
\asmm{n, l, k =0}{\infty} 
 \z{\sfrq{b,c}{n} \sfrq{\z{qa}{bc}}{N-n}\lrdl{(}{)}
{\z{qa}{bc}}^{n}}
{\sfrq{\z{qa}{b}, \z{qa}{c}}{N} (q; q)_{N-n}}
\z{\sfrq{b',c'}{l} \sfrq{\z{qa}{b'\ c'}}{L-l}\lrdl{(}{)}
{\z{qa}{b'\ c'}}^{l} }
{\sfrq{\z{qa}{b'}, \z{qa}{c'}}{L} (q; q)_{L-l}} \\ \\
\times\ \z{\sfrq{b'',c''}{k} \sfrq{\z{qa}{b''\ c''}}{K-k}\lrdl{(}{)}
{\z{qa}{b''\ c''}}^{k}\ \btt_{(n, l, k)}}
{\sfrq{\z{qa}{b''}, \z{qa}{c''}}{K} (q; q)_{K-k}} \\ \\
=\asmm{m=0}{min. (N, L, K)}
\z{\sfrq{b, c, b',\ c',\ b'',\ c'',}{m} 
\lrdl{(}{)}{\z{a^{3}q^{3}}
{b\ c\ b'\ c'\ b''\ c''}}^{m} }
{\sfrq{\z{qa}{b}, \z{qa}{c}, \z{qa}{b'}, \z{qa}{c'}
\z{qa}{b''}, \z{qa}{c''}}{m}
}\\ \\
\times\ \z{\aph_{m}}
{(q; q)_{N-m}(q; q)_{L-m}(q;q)_{K-m}(aq;q)_{N+m}
(aq;q)_{L+m}(aq;q)_{K+m}} . 
\tag{5.5} 
\end{multline} 
Now taking $ b, c, b',\ c',\ b'',\ c'',\  M, N \rightarrow \infty $ in (4.5), 
we get 
$$
  \asmm{n, l, k =0}{\infty} a^{n+l+k}\ 
q^{n^{2}+l^{2}+k^{2}}\ \btt_{(n, l, k)} = \z{1}{\sfrq{aq}{\infty}^{3}}
\ \ \asmm{m=0}{\infty} a^{3m}\ q^{3m^{2}}\ \aph_{m} ,
 \eqno(5.6)$$
for any SBTP $( \aph_{m}, \btt_{(n, l, k)})$. 
Now, we shall make use of the following two ``SBTP" deduced by us : 
\begin{align*} \begin{split} 
&\aph_{m}= \z{\sfrq{a}{m} (1-aq^{2m})}{\sfrq{q}{m} (1-a)} 
(-1)^{m} a^{m}\ q^{\frac{1}{2} (3m^{2}-m)} ,  \\[9pt] 
&\btt_{(n, l, k)}= \z{\sfrq{aq}{n+l+k}}{\sfrq{qa}{n+l} 
\sfrq{qa}{n+k}\sfrq{qa}{l+k}\sfrq{q}{n} \sfrq{q}{l}\sfrq{q}{k}} . 
\end{split} \tag{5.7} 
\end{align*} 
The fact that ``SBTP" $( \aph_{m}, \btt_{(n, l, k)})$ given by (5.7) 
satisfy (5.1) may be verified by substituting it in (5.1) 
and then appealing to \sxpf . \\[9pt]   
\subsection{New triple series \rr\ type \id\ and 
corresponding infinite families  :} 
To derive \rr\ type \id , 
we insert the  ``SBTP" (5.7)  in (5.6) to get (5.8) as given below: 
\begin{multline}
\asmm{n, l, k =0}{\infty} \z{ a^{n+l+k}\ 
q^{n^{2}+l^{2}+k^{2}}\sfrq{qa}{n+l+k}}{\sfrq{qa}{n+l} 
\sfrq{qa}{n+k}\sfrq{qa}{l+k}\sfrq{q}{n} \sfrq{q}{l}\sfrq{q}{k}}\\[9pt]
= \z{1}{\sfrq{aq}{\infty}^{3}}
\asmm{m=0}{\infty} \z{\sfrq{a}{m} (1-aq^{2m})}{\sfrq{q}{m} (1-a)} 
 (-1)^{m} a^{4m}\ q^{\frac{1}{2} (9m^{2}-m)} 
\tag{5.8}
\end{multline}
and then setting $ a=1 $ and $ a=q $ in (5.8)  and
 using Jacobi triple product identity, we obtain following two 
 new triple series \rr\ type \id :
\begin{multline}
\asmm{n, l, k =0}{\infty} \z{  
q^{n^{2}+l^{2}+k^{2}}\sfrq{q}{n+l+k}}{\sfrq{q}{n+l} 
\sfrq{q}{n+k}\sfrq{q}{l+k}\sfrq{q}{n} \sfrq{q}{l}\sfrq{q}{l}} \\ \\ 
= \z{1}{\sfrq{q}{\infty}^{2}} \spdt{0, \pm 4 \text{(mod\ 9)}} ,  
\tag{5.9}
\end{multline} 
\begin{multline}
\asmm{n, l, k =0}{\infty} \z{  
q^{n^{2}+l^{2}+k^{2}+n+l+k}\sfrq{q^{2}}{n+l+k}}
{\sfrq{q^{2}}{n+l}\sfrq{q^{2}}{n+k}\sfrq{q^{2}}{l+k} 
\sfrq{q}{n} \sfrq{q}{l}\sfrq{q}{k}} \\ \\ 
= \z{1}{\sfrq{q^{2}}{\infty}^{2}} \spdt{0, \pm 1 \text{(mod\ 9)}} . 
\tag{5.10}
\end{multline}
\\ 
Further, like classical Bailey lemma and FBTL, the idea of repeated 
application of  SBTL may be expressed in the concept
 of the ``second Bailey type chain"  a sequence of SBTPs:
$$
(\aph_{m}, \btt_{(n, l, k)})\rightarrow (\aph'_{m}, \btt'_{(n, l, k)})\rightarrow 
(\aph''_{m}, \btt''_{(n, l, k)})\rightarrow (\aph'''_{m}, \btt'''_{(n, l, k)})\rightarrow 
\cdots ;
$$
 Using this chain we can establish:
\subsubsection*{Theorem 5.2}
\begin{align*} \hspace{-1cm}
\asmm{m \geq 0}{} & \z{\sfrq{b_{1}}{m}\sfrq{c_{1}}{m} 
\sfrq{b_{2}}{m}\sfrq{c_{2}}{m} \cdots 
\sfrq{b_{s}}{m}\sfrq{c_{s}}{m}}
{\sfrq{\z{aq}{b_{1}}}{m}\sfrq{\z{aq}{c_{1}}}{m} 
\sfrq{\z{aq}{b_{2}}}{m}\sfrq{\z{aq}{c_{2}}}{m} \cdots 
\sfrq{\z{aq}{b_{s}}}{m}\sfrq{\z{aq}{c_{s}}}{m}} 
\\ \\
& \times\ \z{\sfrq{b'_{1}}{m}\sfrq{c'_{1}}{m} 
\sfrq{b'_{2}}{m}\sfrq{c'_{2}}{m} \cdots 
\sfrq{b'_{s}}{m}\sfrq{c'_{s}}{m}}
{\sfrq{\z{aq}{b'_{1}}}{m}\sfrq{\z{aq}{c'_{1}}}{m} 
\sfrq{\z{aq}{b'_{2}}}{m}\sfrq{\z{aq}{c'_{2}}}{m} \cdots 
\sfrq{\z{aq}{b'_{s}}}{m}\sfrq{\z{aq}{c'_{s}}}{m}} \\ \\ 
& \times\ \z{\sfrq{b''_{1}}{m}\sfrq{c''_{1}}{m} 
\sfrq{b''_{2}}{m}\sfrq{c''_{2}}{m} \cdots 
\sfrq{b''_{s}}{m}\sfrq{c''_{s}}{m}}
{\sfrq{\z{aq}{b''_{1}}}{m}\sfrq{\z{aq}{c''_{1}}}{m} 
\sfrq{\z{aq}{b''_{2}}}{m}\sfrq{\z{aq}{c''_{2}}}{m} \cdots 
\sfrq{\z{aq}{b''_{s}}}{m}\sfrq{\z{aq}{c''_{s}}}{m}} \\ \\ 
& \times\ \z{\sfrq{q^{-N}}{m}\sfrq{q^{-L}}{m}\sfrq{q^{-K}}{m}}
{\sfrq{aq^{1+N}}{m}\sfrq{aq^{1+L}}{m}\sfrq{aq^{1+K}}{m}} \\ \\
& \lrdl{(}{)}{\frac{a^{3s}\ q^{3s+N+L+K}}
{b_{1}c_{1} \cdots b_{s}c_{s} b'_{1}c'_{1} 
\cdots b'_{s}c'_{s} b''_{1}c''_{1} 
\cdots b''_{s}c''_{s}}}^{m}
 q^{\frac{3}{2}(-m^{2}+m)}\ \aph_{m} \\ \\ 
& =\z{\sfrq{aq}{N} \sfrq{\z{aq}{b_{s}c_{s}}}{N}}
{\sfrq{\z{aq}{b_{s}}}{N}\sfrq{\z{aq}{c_{s}}}{N}} 
\z{\sfrq{aq}{L} \sfrq{\z{aq}{b'_{s}c'_{s}}}{L}}
{\sfrq{\z{aq}{b'_{s}}}{L}\sfrq{\z{aq}{c'_{s}}}{L}}
\z{\sfrq{aq}{K} \sfrq{\z{aq}{b''_{s}c''_{s}}}{K}}
{\sfrq{\z{aq}{b''_{s}}}{K}\sfrq{\z{aq}{c''_{s}}}{K}}
\\ \\ & \times\ 
\asmm{\substack{n_{s} \geq \cdots \geq n_{1} \geq 0 
\\ \\  l_{s} \geq \cdots \geq l_{1} \geq 0
\\ \\  k_{s} \geq \cdots \geq k_{1} \geq 0
}}{} 
\z{\sfrq{b_{s}}{n_{s}}\sfrq{c_{s}}{n_{s}} 
 \cdots 
\sfrq{b_{1}}{n_{1}}\sfrq{c_{1}}{n_{1}}} 
{\sfrq{q}{n_{s}-n_{s-1}}\sfrq{q}{n_{s-1}-n_{s-2}} 
\cdots \sfrq{q}{n_{2}-n_{1}}}
\end{align*} 
\begin{align*} & \times\
\z{\sfrq{b'_{s}}{l_{s}}\sfrq{c'_{s}}{l_{s}} 
 \cdots 
\sfrq{b'_{1}}{l_{1}}\sfrq{c'_{1}}{l_{1}}}
{\sfrq{q}{l_{s}-l_{s-1}}\sfrq{q}{l_{s-1}-l_{s-2}} 
\cdots \sfrq{q}{l_{2}-l_{1}}}
\\ \\ & \times\
\z{\sfrq{b''_{s}}{k_{s}}\sfrq{c''_{s}}{k_{s}} 
 \cdots 
\sfrq{b''_{1}}{l_{1}}\sfrq{c''_{1}}{l_{1}}}
{\sfrq{q}{k_{s}-k_{s-1}}\sfrq{q}{k_{s-1}-k_{s-2}} 
\cdots \sfrq{q}{k_{2}-k_{1}}}
\\ \\
& \times\ \z{\sfrq{q^{-N}}{n_{s}}}
{\sfrq{\z{b_{s}\ c_{s}q^{-N}}{a}}{n_{s}}}
\z{\sfrq{q^{-L}}{l_{s}}}
{\sfrq{\z{b'_{s}\ c'_{s}q^{-L}}{a}}{l_{s}}}
\z{\sfrq{q^{-K}}{k_{s}}}
{\sfrq{\z{b''_{s}\ c''_{s}q^{-K}}{a}}{k_{s}}}
\\ \\ & \times\ 
\z{ \sfrq{\z{aq}{b_{s-1}\ c_{s-1}}}{n_{s}-n_{s-1}}
\cdots \sfrq{\z{aq}{b_{1}\ c_{1}}}{n_{2}-n_{1}}
}
{\sfrq{\z{aq}{b_{s-1}}}{n_{s}}\sfrq{\z{aq}{c_{s-1}}}{n_{s}} 
\cdots \sfrq{\z{aq}{b_{1}}}{n_{2}}\sfrq{\z{aq}{c_{1}}}{n_{2}}
} 
\\ \\ & \times\ 
\z{ \sfrq{\z{aq}{b'_{s-1}\ c'_{s-1}}}{l_{s}-l_{s-1}}
\cdots \sfrq{\z{aq}{b'_{1}\ c'_{1}}}{l_{2}-l_{1}}
}
{\sfrq{\z{aq}{b'_{s-1}}}{l_{s}}\sfrq{\z{aq}{c'_{s-1}}}{l_{s}} 
\cdots \sfrq{\z{aq}{b'_{1}}}{l_{2}}\sfrq{\z{aq}{c'_{1}}}{l_{2}}
} \\ \\ 
& \times\ 
\z{ \sfrq{\z{aq}{b''_{s-1}\ c''_{s-1}}}{k_{s}-k_{s-1}}
\cdots \sfrq{\z{aq}{b''_{1}\ c''_{1}}}{k_{2}-k_{1}}
}
{\sfrq{\z{aq}{b''_{s-1}}}{k_{s}}\sfrq{\z{aq}{c''_{s-1}}}{k_{s}} 
\cdots \sfrq{\z{aq}{b''_{1}}}{k_{2}}\sfrq{\z{aq}{c''_{1}}}{k_{2}}
} \\ \\ 
& \times\ 
q^{n_{1}+n_{2}+ \cdots +n_{s}+l_{1}+l_{2}+ \cdots +l_{s}
+k_{1}+k_{2}+ \cdots +k_{s}
} \\ \\  
& \times\ a^{n_{1}+n_{2}+ \cdots +n_{s-1}+l_{1}+l_{2}+ \cdots +l_{s-1}
+k_{1}+k_{2}+ \cdots +k_{s-1}} 
\\ \\ & \times\ 
\lrdl{(}{)}{b_{s-1}c_{s-1}}^{-n_{s-1}} \cdots 
\lrdl{(}{)}{b_{1}c_{1}}^{-n_{1}} 
\lrdl{(}{)}{b'_{s-1}c'_{s-1}}^{-l_{s-1}} \cdots 
\lrdl{(}{)}{b'_{1}c'_{1}}^{-l_{1}} \\ \\
& \times\ \lrdl{(}{)}{b'_{s-1}c'_{s-1}}^{-k_{s-1}} \cdots 
\lrdl{(}{)}{b'_{1}c'_{1}}^{-k_{1}}\ 
 \btt_{(n_{1}, l_{1}, k_{1})} . 
\tag{5.11} 
\end{align*}
\subsubsection*{Proof :}
Equation (5.11) is the result of  a $s$-fold iteration of SBTL 
with the choice $ b=b_{i}$, $c=c_{i}, b'=b'_{i}$, $c'=c'_{i},  
b''=b''_{i}$, $c''=c''_{i}$ at the $i$ th step. 
$ \square $ 
\subsubsection*{Theorem 5.3}
If $ (\aph_{m}, \btt_{(n, l, k)}) $ is a SBTP (i.e. related by (5.1)),
 then 
\begin{align*} \hspace{-1cm} 
\asmm{\substack{n_{s} \geq \cdots \geq n_{1} \geq 0 
\\ \\  l_{s} \geq \cdots \geq l_{1} \geq 0
\\ \\  k_{s} \geq \cdots \geq k_{1} \geq 0}}{}  &   
\z{  
q^{n^{2}_{1}+n^{2}_{2}+ \cdots +n^{2}_{s}
+l^{2}_{1}+l^{2}_{2}+ \cdots +l^{2}_{s}
+k^{2}_{1}+k^{2}_{2}+ \cdots +k^{2}_{s}}
}
{\sfrq{q}{n_{s}-n_{s-1}} 
\cdots \sfrq{q}{n_{2}-n_{1}}
\sfrq{q}{l_{s}-l_{s-1}} 
\cdots \sfrq{q}{l_{2}-l_{1}}}\\ \\ 
 & \times\ 
\z{a^{n_{1}+n_{2}+ \cdots +n_{s}+l_{1}+l_{2}+ \cdots +l_{s}
+k_{1}+k_{2}+ \cdots +k_{s}}
\ \btt_{(n_{1}, l_{1}, k_{1})}}
{\sfrq{q}{k_{s}-k_{s-1}} 
\cdots \sfrq{q}{k_{2}-k_{1}}}
\\ \\ 
 &
= \z{1}{\sfrq{aq}{\infty}^{3}}\ \ 
\asmm{m=0}{\infty} q^{3sm^{2}}\ a^{3sm}\ \aph_{m} 
\tag{5.12}
\end{align*} 
\subsubsection*{Proof :}
Letting $ N, L, K,  b_{1}, \cdots, b_{s}, c_{1},  
\cdots c_{s}, b'_{1}, \cdots, b'_{s}, c'_{1},  
\cdots c'_{s}, b''_{1}, \cdots, b''_{s}, c''_{1},  
\cdots c''_{s} $ 
all tend to infinity in Theorem 5.2, the proof of Theorem 5.3 follows. $ \square$  \\[9pt]  
Theorem 5.3 may now be used to embed the triple series \rr\ 
type \id\ (5.9) and (5.10)  in an infinite family of such identities:
\begin{multline*} \hspace{-1.3cm} 
\asmm{\substack{n_{s} \geq \cdots \geq n_{1} \geq 0 
\\ \\  l_{s} \geq \cdots \geq l_{1} \geq 0 
\\ \\  k_{s} \geq \cdots \geq k_{1} \geq 0}}{}  
\z{ 
q^{n^{2}_{1}+ \cdots +n^{2}_{s}
+l^{2}_{1}+ \cdots +l^{2}_{s}+k^{2}_{1}+ \cdots +k^{2}_{s}}
}
{\sfrq{q}{n_{s}-n_{s-1}} 
\cdots \sfrq{q}{n_{2}-n_{1}}
\sfrq{q}{l_{s}-l_{s-1}} 
\cdots \sfrq{q}{l_{2}-l_{1}} }\\ \\
\times\ \z{\sfrq{q}{n_{1}+l_{1}+k_{1}}}
{\sfrq{q}{k_{s}-k_{s-1}} 
\cdots \sfrq{q}{k_{2}-k_{1}}
\sfrq{q}{n_{1}+l_{1}}\sfrq{q}{n_{1}+k_{1}}\sfrq{q}{l_{1}+k_{1}}
} \\ \\ 
\times\ \z{1}{\sfrq{q}{n_{1}}\sfrq{q}{l_{1}}\sfrq{q}{k_{1}}}
= \z{1}{\sfrq{q}{\infty}^{2}}
 \spdt{0, \pm (3s+1) \text{(mod\ 6s+3)}}
\tag{5.13}
\end{multline*} \\  and 
\begin{multline*} \hspace{-1.3cm} 
\asmm{\substack{n_{s} \geq \cdots \geq n_{1} \geq 0 
\\ \\  l_{s} \geq \cdots \geq l_{1} \geq 0 
\\ \\  k_{s} \geq \cdots \geq k_{1} \geq 0}}{}  
\z{ 
q^{n^{2}_{1}+ \cdots +n^{2}_{s}
+l^{2}_{1}+ \cdots +l^{2}_{s}+k^{2}_{1}+ \cdots +k^{2}_{s}}
}
{\sfrq{q}{n_{s}-n_{s-1}} 
\cdots \sfrq{q}{n_{2}-n_{1}}
\sfrq{q}{l_{s}-l_{s-1}} 
\cdots \sfrq{q}{l_{2}-l_{1}} }\\ \\
\times\ \z{q^{n_{1}+ \cdots +n_{s}+l_{1}+ \cdots +l_{s}+
k_{1}+ \cdots +k_{s}}}
{\sfrq{q}{k_{s}-k_{s-1}} 
\cdots \sfrq{q}{k_{2}-k_{1}}
\sfrq{q^{2}}{n_{1}+l_{1}}\sfrq{q^{2}}{n_{1}+k_{1}}\sfrq{q^{2}}{l_{1}+k_{1}}
} \\ \\ 
\times\ \z{\sfrq{q^{2}}{n_{1}+l_{1}+k_{1}}}
{\sfrq{q^{2}}{n_{1}}\sfrq{q^{2}}{l_{1}}\sfrq{q^{2}}{k_{1}}}
= \z{1}{\sfrq{q^{2}}{\infty}^{2}}
 \spdt{0, \pm 1 \text{(mod\ 6s+3)}} . 
\tag{5.14}
\end{multline*} 

\begin{flushleft}
106, Arihant Nagar, Kalka Mata Road,\\
Udaipur-313 001, Rajasthan, INDIA.\\
--- ---\\
Department of Mathematics and Statistics,\\
University College of Science,\\
Mohanlal Sukhadia University,\\
Udaipur-313 001, Rajasthan, INDIA.\\
 E-Mail : {\it yashoverdhan@yahoo.com} \\
\end{flushleft}
\end{document}